\documentclass[12pt]{amsart}
\usepackage{amssymb,amsmath,amsfonts,latexsym,setspace}
\usepackage{bm}
\newcommand{\newsection}[1]{\setcounter{equation}{0} \section{#1}}
\newcommand{\bea}{\begin{eqnarray}}
\newcommand{\eea}{\end{eqnarray}}

\newcommand{\vp}{\varphi}
\newcommand{\cd}{\cdot}

\newcommand{\clb}{\mathcal{B}}

\newcommand{\cle}{\mathcal{E}}
\newcommand{\clf}{\mathcal{F}}

\newcommand{\clh}{\mathcal{H}}
\newcommand{\clk}{\mathcal{K}}

\newcommand{\clm}{\mathcal{M}}

\newcommand{\raro}{\rightarrow}

\def \qed {\hfill \vrule height6pt width 6pt depth 0pt}
\def\textmatrix#1&#2\\#3&#4\\{\bigl({#1 \atop #3}\ {#2 \atop #4}\bigr)}
\def\dispmatrix#1&#2\\#3&#4\\{\left({#1 \atop #3}\ {#2 \atop #4}\right)}
\newcommand{\be}{\begin{equation}}
\newcommand{\ee}{\end{equation}}
\newcommand{\ben}{\begin{eqnarray*}}
\newcommand{\een}{\end{eqnarray*}}

\newcommand{\NI}{\noindent}

\newcommand{\bi}{\begin{itemize}}
\newcommand{\ei}{\end{itemize}}

\newtheorem{Theorem}{\sc Theorem}[section]
\newtheorem{Lemma}[Theorem]{\sc Lemma}
\newtheorem{Proposition}[Theorem]{\sc Proposition}
\newtheorem{Corollary}[Theorem]{\sc Corollary}
\newtheorem{Definition}[Theorem]{\sc Definition}
\newtheorem{Example}[Theorem]{\sc Example}
\newtheorem{Remark}[Theorem]{\sc Remark}
\newtheorem{Note}[Theorem]{\sc Note}
\newtheorem{Question}{\sc Question}
\newtheorem{ass}[Theorem]{\sc Assumption}
\newcommand{\bt}{\begin{Theorem}}
\def\beginlem{\begin{Lemma}}
\def\beginprop{\begin{Proposition}}
\def\begincor{\begin{Corollary}}
\def\begindef{\begin{Definition}}
\def\beginexamp{\begin{Example}}
\def\beginrem{\begin{Remark}}
\def\beginq{\begin{Question}}
\def\beginass{\begin{ass}}
\def\beginnote{\begin{Note}}
\newcommand{\et}{\end{Theorem}}
\def\endlem{\end{Lemma}}
\def\endprop{\end{Proposition}}
\def\endcor{\end{Corollary}}
\def\enddef{\end{Definition}}
\def\endexamp{\end{Example}}
\def\endrem{\end{Remark}}
\def\endq{\end{Question}}
\def\endass{\end{ass}}
\def\endnote{\end{Note}}

\begin{document}

\title[Curvature invariant and generalized canonical operator models - II]{Curvature invariant and generalized canonical operator models - II}
\author[Douglas]{Ronald G. Douglas}
\address[Ronald G. Douglas]{Department of Mathematics, Texas A \& M University, College Station, TX, 77843, USA}
\email{rdouglas@math.tamu.edu}

\author[Kim]{Yun-Su Kim}
\address[Yun-Su Kim]{Deceased}

\author[Kwon]{Hyun-Kyoung Kwon}
\address[Hyun-Kyoung Kwon]{Department of Mathematics, The University of Alabama, Tuscaloosa, AL, 35487, USA}
\email{hkwon5@as.ua.edu}

\author[Sarkar]{Jaydeb Sarkar}
\address[Jaydeb Sarkar]{Statistics and Mathematics Unit, Indian Statistical Institute, Bangalore, 560059, India}
\email{jay@isibang.ac.in, jaydeb@gmail.com}

\thanks{The work of Douglas and Sarkar was partially supported by a grant from the National
Science Foundation. The work of Kwon was supported in part by the
College Academy of Research, Scholarship and Creative Activity of
the University of Alabama, the Basic Science Research Progra through
the National Research Foundation of Korea (NRF) funded by the
Ministry of Education, Science, and Technology (2011-0026989), the
TJ Park Postdoctoral Fellowship, and by a Young Investigator Award
at the NSF-sponsored Workshop in Analysis and Probability, Texas A
\& M University, 2009. Sarkar would also like to acknowledge the
hospitality of the Department of Mathematics at Texas A \& M
University and at the University of Texas at San Antonio, where part
of his research was done.}



\begin{abstract}
In \cite{DKKS}, the authors investigated a family of quotient Hilbert
modules in the Cowen-Douglas class over the unit disk constructed
from classical Hilbert modules such as the Hardy and Bergman
modules. In this paper we extend the results to the multivariable
case of higher multiplicity. Moreover, similarity as well as
isomorphism results are obtained.
\end{abstract}

\dedicatory{\textsf{Dedicated to the memory of our colleague Yun-Su
Kim}}

\maketitle

\newsection{Introduction}
In \cite{DKKS}, we reported on some results obtained in comparing
the canonical model of Sz.-Nagy and Foias with the complex geometric
model of M. J. Cowen and the first author. These models are first
recast in the language of Hilbert modules as in \cite{DP}. In
\cite{DKKS}, we considered only the simplest non-trivial cases of
quotient modules, first of the Hardy module mapped into the
$\mathbb{C}^2$-valued Hardy module. Afterwards we extended the
results by replacing the Hardy module by other Hilbert modules over
$\mathbb{C}[z]$ related to the unit disk $\mathbb{D}$ such as the
Bergman and weighted Bergman modules. We also took up the question
of when two such quotient modules are isomorphic.

In this paper we proceed to the more general cases of these
phenomena in higher multiplicity and for Hilbert modules over
$\mathbb{C}[z_1,\cdots, z_n]$ and again, determine when two such
quotient Hilbert modules are (isometrically) isomorphic and, in some
cases, similar. Here we represent the associated hermitian
anti-holomorphic vector bundle as a twisted tensor product of the
vector bundle for the basic Hilbert module by a vector bundle
determined by the multiplier used to define the quotient. A version
of this representation was used earlier by Uchiyama \cite{U} and
Treil and the third author \cite{KT}. (See \cite{Zhu} also for some
related results concerning Hilbert modules over $\mathbb{C}[z]$.) We
observe that although the vector bundles in the short exact sequence
of vector bundles defining the quotient module are all pull-backs
from an infinite dimensional Grassmanian and hence, appear to be
built using an infinite dimensional Hilbert space, they are all
actually the tensor product of a quotient of finite rank, trivial
vector bundles by the fixed line bundle for the basic Hilbert
module. Thus, all calculations and proofs can be carried out in this
finite dimensional context.

After some preliminaries in Section 2, we describe the results in
\cite{DKKS} and take up the issue of isomorphism for a more general
case of multiplicity one. In Section 3, we extend the definitions of
Section 2 to the multivariate case and then proceed to extend and
generalize the key results from Section 2 to cases of higher
multiplicity.

In Section 4, we explore some similarity questions for quotient
Hilbert modules in this context drawing on the research of two
earlier groups. First, we use results on the similarity question
in the Hardy space context which originated in the research of
Sz.-Nagy and Foias \cite{NF} and, more recently, of Treil and the
third author \cite{KT}. In the latter work, similarity is shown to
be equivalent to the existence of a certain bounded function whose
Laplacian is related to the curvature of the quotient vector bundle.
Second, a study by a group of Chinese researchers (cf. \cite{JW})
showed that in the case of contractive Hilbert modules over
$\mathbb{C}[z]$, some results for similarity are independent of the
particular basic Hilbert module used to construct the quotient
Hilbert module. For example, a quotient Hilbert module defined using
the Bergman module is similar to the Bergman module itself if and
only if the same is true for the analogous quotient Hilbert module
defined using the Hardy module. Our proof of this fact rests on the
tensor product factorization mentioned above since the finite
dimensional vector bundles involved do not depend on the basic
Hilbert module used. A related result based on techniques from
function theory appears in \cite{DKT}. Finally in Section 5, we
raise some questions suggested by these results and point out
further connections with other work.

\newsection{Preliminaries}

We first consider several definitions and facts concerning Hilbert
modules. We denote by  $\mathbb{C}[z_1, \ldots, z_n]$ the algebra of
polynomials in $n$ commuting variables $z_1, \ldots, z_n$.

Now let $\clh$ be a Hilbert space and let $\{T_1, \ldots, T_n\}
\subseteq \clb(\clh)$ be an $n$-tuple of commuting operators on
$\clh$. Then the operators $\{T_1, \ldots, T_n\}$ induce a module
action on the Hilbert space $\clh$ over $\mathbb{C}[z_1, \ldots,
z_n]$ as follows (see \cite{DP}):
\[
p \cdot h := p(T_1, \ldots, T_n) h,
\]
for all $p \in \mathbb{C}[z_1, \ldots, z_n]$ and $h \in \clh$.
Denote by $M_p: \clh \raro \clh$ the bounded linear operator
\[
M_p h = p \cdot h = p(T_1, \ldots, T_n)h,
\]
for $h \in \clh$. In particular, for $p = z_i \in \mathbb{C}[z_1
\ldots, z_n]$, $i=1,\ldots, n$, we obtain the module multiplication
operators $M_{z_i}$ defined by
$$M_{z_i} h = z_i \cdot h = T_i h,$$
for $h \in \clh$. In what follows, we will use the notion of a
Hilbert module $\clh$ over $\mathbb{C}[z_1, \ldots, z_n]$ in place
of an $n$-tuple of commuting operators $\{T_1, \ldots, T_n\}
\subseteq \clb(\clh)$, where the operators are determined by module
multiplication by the coordinate functions, and vice versa.

The notion of intertwining maps between operators on Hilbert spaces
can be formulated in terms of module maps between Hilbert modules.

\begin{Definition}
A bounded linear map $X: \clh \rightarrow \tilde{\clh}$ between two
Hilbert modules $\clh$ and $\tilde{\clh}$ over$\text{ }$  $\mathbb{C}[z_1.
\ldots, z_n]$ is said to be a module map if $X M_{z_i} = M_{z_i} X$
for $i = 1, \ldots, n$, or equivalently, if $X M_p = M_p X$ for $p
\in \mathbb{C}[z_1, \ldots, z_n]$. Two Hilbert modules are said to
be isomorphic if there exists a unitary module map between them, and
similar if there exists an invertible module map between them.
\end{Definition}

A natural source of Hilbert modules is the family of reproducing
kernel Hilbert spaces \cite{A} on domains in $\mathbb{C}^n$ with
bounded multiplication operators defined by the coordinate
functions.

Recall that for a non-empty set $X$ and a Hilbert space $\cle$, an
operator-valued function $K: X \times X \rightarrow \clb(\cle)$ is
said to be a \textit{positive definite kernel} if \[ \sum_{i, j =
1}^{k} \langle K(z_i, z_j) \eta_j, \eta_i \rangle > 0,\] for all
$\eta_i \in \cle, \, z_i \in X, i=1, \cdots, k$, and all $k \in
\mathbb{N}$. Given such a positive definite kernel $K$ on $X$, one
can construct a Hilbert space $\clh_K$ of $\cle$-valued functions on
$X$ as the completion of the linear span of the set $\{K(\cdot, w)
\eta : w \in X, \eta \in \cle\}$ with respect to the inner product
\[\langle K(\cdot, w) \eta, K(\cdot, z) \zeta\rangle_{\clh_K} = \langle K(z,
w) \eta, \zeta\rangle_{\cle},\]for all $z, w \in X$ and $\eta, \zeta
\in \cle$. The kernel function $K$ has the reproducing property so
that for $f \in \clh_K$ and $\eta \in \cle$, one has
\[\langle f, K(\cdot, z) \eta\rangle_{\clh_K} = \langle f(z),
\eta\rangle_{\cle}.\]Hence, the evaluation operator $\bm{ev}_z :
\clh_K \rightarrow \cle$ defined by \[\langle \bm{ev}_z(f),
\eta\rangle_{\cle} = \langle f, K(\cdot, z) \eta\rangle_{\clh_K}\]
is bounded for all $z \in X$, where $\eta \in \cle$ and $f \in
\clh_K$.

Conversely, given a Hilbert space $\clh$ of functions from $X$ to
$\cle$ with bounded evaluation operators $\bm{ev}_z$ for each $z \in
X$, one can define a kernel function $K$ corresponding to $\clh$ as
\[K(z, w) = \bm{ev}_z \circ \bm{ev}_w^* \in \clb(\cle),\]for all $z, w \in X$. If
$X$ is a domain $\Omega$ in $\mathbb{C}^n$ and $K : \Omega
\times \Omega \rightarrow \cle$ is holomorphic in the first variable
and anti-holomorphic in the second variable, then $\clh$ is a
space of holomorphic functions on $\Omega$. If, in addition, the
multiplication operators $M_{z_1}, \ldots, M_{z_n}$ by the
coordinate functions are bounded on $\clh$, then we say that
$\clh$ is a \textit{reproducing kernel Hilbert module} over
$\Omega$. In this case, it is easy to verify that \[M_{z_i}^*
(K(\cdot, w) \eta) = \bar{w}_i K(\cdot, w) \eta,\]for all $w = (w_1,
\ldots, w_n) \in \Omega, \eta \in \cle$, and $i=1, \ldots, n$. It is
easy to show that a necessary condition for $\clh$ to be a Hilbert
module is for $\Omega$ to be bounded.

Given $\cle$- and $\cle_*$-valued reproducing kernel Hilbert modules
$\clh$ and $\clh_*$, respectively, over the domain $\Omega$, a
function $\varphi : \Omega \rightarrow \clb(\cle, \cle_*)$ is said
to be a \textit{multiplier} if $\varphi f \in \clh_*$ for all $f \in
\clh$, where $(\varphi f)(z) = \varphi(z) f(z)$ for all $z \in
\Omega$. The set of all such multipliers is denoted $\clm(\clh,
\clh_*)$, or simply $\clm(\clh)$, if $\clh= \clh_*$. By the closed graph theorem, each $\varphi \in \clm(\clh,
\clh_*)$ induces a bounded linear map $M_{\varphi} : \clh
\rightarrow \clh_*$. Consequently, $\clm(\clh, \clh_*)$ is a Banach
space with
\[\|\varphi\|_{\clm(\clh, \clh_*)} = \|M_{\varphi}\|_{\clb(\clh,\clh_*)}.\]
\NI For $\clh = \clh_*$, $\clm(\clh)$ is a Banach
algebra with this norm. One can also view a multiplier as a module
map $\Theta \in \clb(\clh, \clh_*)$ since such an operator is given
by pointwise multiplication by a function from $\Omega$ to
$\clb(\cle, \cle_*)$.

A class of reproducing kernel Hilbert modules over $\Omega \subseteq
\mathbb{C}$, denoted by $B_m(\Omega)$, was
introduced by M. J. Cowen and the first author in \cite{CD}. This
notion was extended to the multivariable setting ($\Omega \subseteq
\mathbb{C}^n$) in \cite{CD2}, by Curto and Salinas \cite{CS}, and by Chen and the
first author \cite{ChenD} to the Hilbert module context. We focus on
the dual $B^*_m(\Omega)$ of $B_m(\Omega)$.

\begin{Definition}\label{CD-Defn}
Let $\Omega$ be a domain in $\mathbb{C}^n$ and let $m$ be a positive
integer. Then $B^*_m(\Omega)$ is the set of all Hilbert modules
$\clh$ over $\mathbb{C}[z_1, \ldots, z_n]$ such that

(i) The column operator $(M_z - w I_{\clh})^* : \clh \rightarrow
\clh^n$ defined by

\[(M_z - w I_{\clh})^* h = (M_{z_1} - w_1 I_{\clh})^* h \oplus
\cdots \oplus (M_{z_n} - w_n I_{\clh})^* h,\quad \quad (h \in
\clh)\] has closed range for all $w=(w_1, \cdots, w_n) \in \Omega$, where $\clh^n =
\clh \oplus \cdots \oplus \clh$.

(ii) $\mbox{dim~} \mbox{ker~}(M_z - w I_{\clh})^* = \mbox{dim}
[~\cap_{i=1}^n \mbox{ker~}(M_{z_i} - w_i I_{\clh})^*~] = m$ for all $w
\in \Omega$, and

(iii) $\bigvee_{{w} \in \Omega} \mbox{ker~}(M_z - w I_{\clh})^* =
\clh$.
\end{Definition}

Note that these modules are the duals of those in $B_m(\Omega)$
defined in \cite{CD}. The use of this class results in
anti-holomorphic objects as opposed to holomorphic ones but there is no essential difference.

Given a Hilbert module $\clh$ in $B_m^*(\Omega)$, the mapping $w
\mapsto E^*_{\clh}(w) := \mbox{ker~} (M_z - w I_{\clh})^*$ defines a
rank $m$ hermitian anti-holomorphic vector bundle over $\Omega$
which will be denoted by $E^*_{\clh}$. If $E^*_{\clh}$ is trivial,
there exists a $\clb(\mathbb{C}^m)$-valued kernel function on
$\Omega$. More precisely, let $\{s^i_z : 1 \leq i \leq m\}$ be an
anti-holomorphic frame for the vector bundle $E^*_{\clh}$. A kernel
function for $\clh$ can be obtained as the Gram matrix of the frame
$\{s^i_w : 1 \leq i \leq m \}$; that is,
\[K(z, w) = \big( \langle s^j_w, s^i_z\rangle_{E^*_{\clh}(w)})_{i,j=1}^m.\]

If $E^*_{\clh}$ is not trivial, then we can use an anti-holomorphic
frame over an open subset $U \subseteq \Omega$ to define a kernel
function $K_U$ on $U$. Since a domain is connected, one can show
that $\clh_{K_U} \cong \clh$. One way to obtain a local frame is to
identify the fiber of the dual vector bundle $E_{\clh}$ with
$\clh/I_w \cdot \clh \cong \mathbb{C}^m \cong \mbox{span} \{s^i_w :
1 \leq i \leq m\}$, as hermitian Hilbert modules where $I_w = \{ p
\in \mathbb{C}[z_1, \ldots, z_n], \text{ }p(w) = 0\}$ and $w \in \Omega$.

Although one usually considers the Chern connection for hermitian
holomorphic vector bundles, an analogous definition can be used for
hermitian, anti-holomorphic vector bundles. Once one identifies a
fixed basis in a fiber locally, the difference between the two
notions amounts to taking the complex conjugate.

The \textit{curvature} of the bundle $E^*_{\clh}$ for the Chern connection
determined by the metric defined by the Gram matrix, or, if
$E^*_{\clh}$ is not trivial, then with the inner product on
$E^*_{\clh}(w) = \mbox{ker~}(M_z - w I_{\clh})^* \subseteq \clh$, is
given by
\[\clk_{E^*_{\clh}} (w) = (\bar{\partial}_j \{K(w, w)^{-1} \partial_i K(w, w)\})_{i,j=1}^n,\]
for all $w \in \Omega$. Note that the representation of the
curvature matrix defined above is with respect to the basis of
two-forms $\{dz_i \wedge d\bar{z}_j : 1 \leq i,j \leq n\}$. In
particular, for a line bundle, that is, when $m=1$, the curvature
form is given by
\[\begin{split}\clk_{E^*_{\clh}}(w) & = \bar{\partial} K(w,w)^{-1} \partial K(w, w) \\ & = -  \partial \bar{\partial} \log \|K(\cdot, w)\|^2 \\ & = - \sum_{i, j =
1}^n \frac{\partial}{\partial w_i}\frac{\partial}{\partial\bar{w_j}}
\log K(w,w) dw_i \wedge d\bar{w}_j,\end{split}\]for all $w \in
\Omega$.

Since curvature has a coordinate free meaning, we can use a local
frame for calculations. Moreover, curvature is given by a
self-adjoint matrix. Once we fix the basis for the two-forms and
identify the dual of a fiber with itself in the usual manner, it
follows that the curvature for the Chern connection on the dual
vector bundle $E^*_{\clh}$ is the same as that for $E_{\clh}$.

The Hardy module $H^2(\mathbb{B}^n)$, the Bergman module
$A^2(\mathbb{B}^n)$, the weighted Bergman modules
$A^2_{\alpha}(\mathbb{B}^n)$ over the unit ball $\mathbb{B}^n$ and
the Drury-Arveson module $H^2_n$ are all in $B^*_1(\mathbb{B}^n)$. A
further source of Hilbert modules in $B^*_m(\Omega)$ is a family
of quotient Hilbert modules, where the standard examples are
used as building blocks. In particular, in \cite{DKKS}, the authors
considered certain quotient Hilbert modules of $H^2(\mathbb{D}),
A^2(\mathbb{D})$ and $A^2_{\alpha}(\mathbb{D})$ which are all in
$B^*_1(\mathbb{D})$ and determined when two such quotient Hilbert
modules are similar and isomorphic. More precisely, given a pair of
functions $\{\theta_1, \theta_2\}$ in the algebra
$H^{\infty}(\mathbb{D})$ of bounded holomorphic functions on
$\mathbb{D}$, which satisfies the corona condition
\[|\theta_1(z)|^2+|\theta_2(z)|^2 \geq \epsilon
> 0,\]for all $z \in \mathbb{D}$, one can define $\Theta \in
H^{\infty}(\mathbb{D}) \otimes \clb(\mathbb{C}, \mathbb{C}^2)$ by the row
matrix $\Theta(z) = (\theta_1(z), \theta_2(z))$. Then for $\clh =
H^2(\mathbb{D}), A^2(\mathbb{D}),$ or $A^2_{\alpha}(\mathbb{D})$, the
quotient Hilbert module $\clh_{\Theta}$ given by the short exact
sequence
\[0
\longrightarrow \clh  \otimes \mathbb{C}
\stackrel{M_{\Theta}}\longrightarrow \clh \otimes \mathbb{C}^2
\stackrel{\pi_{\Theta}}\longrightarrow \clh_{\Theta} \longrightarrow
0,\]is in $B^*_1(\mathbb{D})$. Here $M_{\Theta}$ is the operator
defined as $M_{\Theta} h = \theta_1 h \otimes e_1 +\theta_2 h
\otimes e_2$ for $h \in \clh$, where $\{e_1, e_2\}$ is an
orthonormal basis for $\mathbb{C}^2$, and $\pi_{\Theta}$ is the
quotient map. Results concerning similarity and isomorphism of such
quotient modules were obtained where the notion of curvature played
a crucial role.

\newsection{Certain quotient modules in $B_1^*(\Omega)$}

One can extend the results described in the previous section to the
case of the quotient Hilbert module $\clh_{\Theta}$, where $\clh \in
B_1^*(\Omega)$ for $\Omega \subseteq \mathbb{C}^n$ and $\Theta$
takes values in $\clb(\mathbb{C}^l,\mathbb{C}^{l+1})$ for $l \in
\mathbb{N}$. We assume that the function $\Theta$ is a multiplier
for $\clh$, which is equivalent to saying that \[M_{\Theta} M_{z_i}
= M_{z_i} M_{\Theta},\]for all $i = 1, \ldots, n$. Two hermitian
holomorphic vector bundles $E_{\clh}$ and $E_{\tilde{\clh}}$ over
$\Omega$ are said to be \textit{equivalent} if there exists a biholomorphic
bundle map between $E_{\clh}$ and $E_{\tilde{\clh}}$ which defines
an isometric isomorphism between the fibers. We will denote this
notion of isomorphism by $E_{\clh} \cong E_{\tilde{\clh}}$.

We finesse the issue of the corona condition used above by assuming
that $\Theta$ has a left inverse which is also in the multiplier
algebra.

\begin{Theorem}\label{mcase}
Let $\clh \in B_1^*(\Omega)$ for $\Omega \subseteq \mathbb{C}^n$ and
let $\Theta : \Omega \raro \clb(\mathbb{C}^l, \mathbb{C}^{l+1})$ for
some positive integer $l$ be such that $\Theta \in \clm(\clh) \otimes
\clb(\mathbb{C}^l, \mathbb{C}^{l+1})$ with a left inverse
$\Psi \in \clm(\clh) \otimes \clb(\mathbb{C}^{l+1},
\mathbb{C}^l) $. Denote by $\clh_{\Theta}$ the quotient Hilbert
module $\clh_{\Theta} = (\clh \otimes \mathbb{C}^{l+1})/ M_{\Theta}
(\clh \otimes \mathbb{C}^l)$. Then

(1) $\clh_{\Theta} \in B_1^*(\Omega)$,

(2) $L_{\Theta}(w) = \mathbb{C}^{l+1}/  \Theta(w) \mathbb{C}^l$
defines a holomorphic line bundle $L_\Theta=\coprod_{w \in \Omega}
L_\Theta(w)$ such that $E^*_{\clh_{\Theta}} \cong E^*_{\clh} \otimes
L^*_{\Theta}$, where $L^*_{\Theta}$ is the hermitian
anti-holomorphic vector bundle dual to the hermitian holomorphic
vector bundle $L_{\Theta}$, and

(3) $\clk_{E^*_{\clh_{\Theta}}} - \clk_{E^*_{\clh}} =
\clk_{L^*_{\Theta}}$ as two-forms.
\end{Theorem}

\NI \textsf{Proof.} Since the statements are all local, we can
proceed pointwise as follows. The fact that $M_{\Theta}$ is left
invertible implies that $M_{\Theta}$ has closed range. Thus one
obtains the short exact sequence
\[0 \longrightarrow \clh  \otimes \mathbb{C}^l
\stackrel{M_{\Theta}}\longrightarrow \clh \otimes \mathbb{C}^{l+1}
\stackrel{\pi_{\Theta}}\longrightarrow \clh_{\Theta} \longrightarrow
0,\] where \[M_{\Theta}
\begin{bmatrix} h_1\\ \vdots\\h_l \end{bmatrix} = [\theta_{i,j}]_{(l+1) \times l}  \begin{bmatrix} h_1\\ \vdots\\h_l \end{bmatrix}
= \begin{bmatrix} \sum_{j=1}^l \theta_{1,j} h_j\\ \vdots\\
\sum_{j=1}^l \theta_{l+1,j} h_j
\end{bmatrix},
\] for $h_i \in \clh$, $i=1, \ldots, l$, and $\pi_{\Theta}$ is the quotient map. Localizing this module sequence at $w \in
\Omega$, that is, taking quotients by $I_w \cdot (\clh \otimes
\mathbb{C}^l)$,  $I_w \cdot (\clh \otimes \mathbb{C}^{l+1}),$ and
$I_w \cdot \clh_{\Theta}$, respectively, we have that the sequence
$$\mathbb{C}_w \otimes \mathbb{C}^l \stackrel{I_{\mathbb{C}_w} \otimes \Theta(w)} \longrightarrow \mathbb{C}_w \otimes
\mathbb{C}^{l+1} \stackrel{\pi_{\Theta}(w)} \longrightarrow
\clh_{\Theta}/ I_w \cdot \clh_{\Theta} \longrightarrow 0$$is exact
\cite{DP}. Here, $\mathbb{C}_w$ is the Hilbert module with module
multiplication defined by $p \cdot \lambda = p(w) \lambda$ for all
$p \in \mathbb{C}[z_1, \ldots, z_n]$ and $\lambda \in \mathbb{C}$.
Since $\mbox{dim } \mbox{ran } \Theta(w)=l$ for all $w \in \Omega$,
it follows that $\mbox{dim } \mbox{ker } \pi_{\Theta}(w) = l$, and
thus
\[\mbox{dim~} \clh_{\Theta}/ I_w \cdot \clh_{\Theta} = \mbox{dim~}
\clh_{\Theta}\Big/\Big(\sum_{i=1}^n (M_{z_i} - w_i I_{\clh})
\clh_{\Theta}\Big) =1,\] for all $w \in \Omega$. Consequently,
\[\mbox{dim~} [\mathop{\cap}_{i=1}^n \mbox{ker~} (M_{z_i} - w_i
I_{\clh})^*|_{\clh_{\Theta}}] = 1.\]

\NI To show that $\clh_{\Theta} \in B_1^*(\Omega),$ we must also
demonstrate that \[\mathop{\bigvee}_{w \in \Omega}
[\mathop{\cap}_{i=1}^n \mbox{ker~} (M_{z_i} - w_i
I_{\clh})^*|_{\clh_{\Theta}}] = \clh_{\Theta}.\]

\NI To this end, let $\{e_i\}_{i=1}^{l+1}$ be the standard
orthonormal basis for $\mathbb{C}^{l+1}$ and let $\Delta_{\Theta}$ be
the formal determinant
\[
\Delta_{\Theta}(w) = \mbox{det}\
\begin{bmatrix} e_1 & \theta_{1,1}(w) & \cdots & \theta_{1,l}(w) \\ \vdots & \vdots & \vdots & \vdots
\\ e_{l+1} & \theta_{l+1,1}(w) & \cdots & \theta_{l+1,l}(w) \end{bmatrix} \in \mathbb{C}^{l+1},
\]
for $w \in \Omega$. Since $\Theta(w)$ has a left inverse
$\Psi(w)$, it follows that $\mbox{rank~} \Theta(w)= l$, and hence
$\Delta_{\Theta}(w) \neq 0$ for all $w \in \Omega$. Set $\gamma_w :=
f_w \otimes \overline{\Delta_{\Theta}(w)} \neq 0$ for all $w \in
\Omega$, where $f_w$ is any non-zero vector in $E^*_{\clh}(w)
\subseteq \clh$ and $\overline{\Delta_{\Theta}(w)}$ is the complex
conjugate of $\Delta_{\Theta}(w)$ relative to the basis
$\{e_i\}_{i=1}^{l+1}$. (Note that since $\mbox{dim~}E^*_{\clh}(w) =
1$, $\gamma_w$ is well-defined up to a non-zero scalar.) Moreover,
consider the inner product of $\gamma_w$ with
\[M_{\Theta} \begin{bmatrix}h_1\\ \vdots\\h_l\end{bmatrix} =
\begin{bmatrix}\sum_{j=1}^l \theta_{1,j} h_j\\\vdots\\\sum_{j=1}^l
\theta_{l+1, j} h_j\end{bmatrix} \in \clh \otimes \mathbb{C}^{l+1},
\]for $\{h_i\}_{i=1}^l \subseteq \clh$. Evaluating the resulting
functions at $w \in \Omega$, we get that these functions are the sum
of the products of $h_i(w)$ with coefficients equal to the
determinants of matrices with repeated columns and hence

\[\langle M_{\Theta} \begin{bmatrix}h_1\\ \vdots\\h_l\end{bmatrix}, \gamma_w \rangle = 0. \] Thus, $\gamma_w \perp \mbox{ran~}M_\Theta$ for all $w \in \Omega$.
Also, it is easy to see that
\[(M_{z_i}^* \otimes I_{\mathbb{C}^{l+1}}) \gamma_w = \bar{w_i}
\gamma_w, \]for $w \in \Omega$ and for all $i = 1, \ldots, n$, so
that
\[\mathop{\cap}_{i=1}^n \mbox{ker~} (M_{z_i} \otimes I_{\mathbb{C}^{l+1}} - w_i I_{\clh \otimes \mathbb{C}^{l+1}})^*|_{\clh_{\Theta}} = \mathbb{C}\cdot
\gamma_w,\]for all $w \in \Omega$.

Now we prove that $\bigvee_{w \in \Omega}f_w \otimes
\overline{\Delta_{\Theta}(w)}= {\clh_{\Theta}}$. (Note that this
space is independent of the particular $f_w$'s chosen.) For all $g
=\sum_{i=1}^{l+1} g_i \otimes e_i \in \clh \otimes
\mathbb{C}^{l+1}$ with $g \perp \gamma_w$ for every $w \in \Omega$,
we must exhibit the representation $g_i(w) = \sum_{j=1} ^l
\eta_j(w) \theta_{ij}(w)$ for $i=1,...,l+1$, where the
$\{\eta_j\}_{j=1}^l$ are functions in $\clh$. Fix $w_0 \in \Omega$.
The assumption $\langle g, \gamma_{w_0} \rangle =0$ implies that
\begin{equation}\label{3.1}
\mbox{det}\ \begin{bmatrix} g_1(w_0) & \theta_{1,1}(w_0) & \cdots
& \theta_{1,l}(w_0) \\ \vdots & \vdots & \vdots & \vdots \\
g_{l+1}(w_0) & \theta_{l+1,1}(w_0) & \cdots & \theta_{l+1,l}(w_0)
\end{bmatrix}=0.
\end{equation}
Now view the matrix
$$\Theta(w_0)=\begin{bmatrix} \theta_{1,1}(w_0) & \cdots &
\theta_{1,l}(w_0) \\ \vdots & \vdots & \vdots \\
\theta_{l+1,1}(w_0) & \cdots & \theta_{l+1,l}(w_0) \end{bmatrix}$$
as the coefficient matrix of a linear system of $(l+1)$ equations in
$l$ unknowns. Since $\mbox{rank~}\Theta(w_0)=l$, some principal
minor (which means taking some $l$ rows) has a non-zero determinant.
Hence, using Cramer's rule, we can uniquely solve for
$\{\eta_j(w_0)\}_{j=1}^l \subseteq \mathbb{C}^l$, at least for these
$l$ rows. But by (\ref{3.1}), the solution must also satisfy the
remaining equation. Hence we obtain the $\{\eta_j(w_0)\}_{j=1}^l
\subseteq \mathbb{C}^l$ and define
$$\Xi(w_0)=\sum_{j=1}^{l} \eta_j(w_0) \otimes e_j,$$
so that
$$
g(w_0)=\Theta(w_0)\Xi(w_0),
$$
for each $w_0 \in \Omega$. After doing this for each $w \in \Omega$,
we use the left inverse $\Psi(w)$ for $\Theta(w)$ to obtain
$$\Xi(w) =(\Psi(w) \Theta(w))\Xi(w)=\Psi(w)(\Theta(w)
\Xi(w))=\Psi(w) g(w) \in \clh \otimes \mathbb{C}^l.$$ Consequently,
$\{\eta_j\}_{j=1}^l \subseteq \clh$ and $\bigvee_{w \in \Omega}
\gamma_w = \clh_{\Theta}$.

Lastly, the closed range property of $\clh_{\Theta}$ follows from
that of $\clh$. In particular, since the column operator $M_z^* -
\bar{w}I_{\clh}$ of Definition \ref{CD-Defn} acting on $\clh \otimes \mathbb{C}^{l+1}$
has closed range and a finite
dimensional kernel, it follows that restricting it to the invariant
subspace $\clh_{\Theta} \subseteq \clh \otimes \mathbb{C}^{l+1}$
yields an operator with closed range.

The proofs of parts (2) and (3) are identical to those of the
analogous statements of Theorem \ref{n-curv} below. \qed

\begin{Corollary}\label{cor-3.2}
Let $\clh \in B_1^*(\Omega)$, where $\Omega \subseteq \mathbb{C}^n$,
 and assume that $\Theta_i: \Omega \raro \clb(\mathbb{C}^l,
\mathbb{C}^{l+1}), \Theta_i \in \clm(\clh) \otimes \clb(\mathbb{C}^l,
\mathbb{C}^{l+1})$, for $i=1,2,$ have left inverses in
$\clm(\clh) \otimes \clb(\mathbb{C}^{l+1}, \mathbb{C}^{l})$.
Then the quotient Hilbert modules $\clh_{\Theta_1}$ and
$\clh_{\Theta_2}$ are isomorphic if and only if \[
\bigtriangledown^2 \mbox{log~} \|\Delta_{\Theta_1}\| =
\bigtriangledown^2 \mbox{log~}\|\Delta_{\Theta_2}\|.\]
\end{Corollary}

\NI\textsf{Proof.} We can choose a $k_w$ so that $k_w \otimes
\overline{\Delta_{\Theta_i}(w)}$, $i=1, 2$, are anti-holomorphic
local cross-sections of $E^*_{\clh_{\Theta_1}}$ and
$E^*_{\clh_{\Theta_2}}$, respectively, over some open subset $U
\subseteq \Omega$. Since every $w_0 \in \Omega$ is contained in such
an open subset $U$ of $\Omega$, we can use (3) of Theorem 3.1 and
the result of \cite{CD} stating that two Hilbert modules
$\clh_{\Theta_1}, \clh_{\Theta_2} \in B^*_1(\Omega)$ are isomorphic if
and only if
$$\clk_{E^*_{\clh_{\Theta_1}}}(z) = \clk_{E^*_{\clh_{\Theta_2}}}(z),$$
for every $z \in \Omega$ to complete the proof. \qed

For finite dimensional spaces $\cle$ and $\cle_*$, and a multiplier
$\Theta: \Omega \raro \clb(\cle, \cle_*)$ with constant rank, one
can define the holomorphic kernel and co-kernel bundles with fibers
$\mbox{ker}\, \Theta(w)$ and $\mbox{coker}\, \Theta(w) = \cle_*/
\Theta(w) \cle$ for $w \in \Omega$, respectively. Moreover, related
Hilbert modules with $\clh \in B^*_m(\Omega)$ can be defined for an
arbitrary $m \geq 1$. Here we consider the simplest case, when $m=1$
and $\mbox{ker~} \Theta(w)=\{0\}$, and obtain only some of the most
direct results.

\begin{Theorem}\label{n-curv}
Let $\clh \in B_1^*(\Omega)$ for $\Omega \subseteq \mathbb{C}^n$ and let
$\Theta: \Omega \raro \clb(\mathbb{C}^p, \mathbb{C}^q),\Theta \in
\clm(\clh) \otimes \clb(\mathbb{C}^p, \mathbb{C}^q)$, have a
left inverse $\Psi \in \clm(\clh) \otimes \clb(\mathbb{C}^q,
\mathbb{C}^p)$, where $1 \leq p < q < \infty$. Then there exists a
hermitian anti-holomorphic vector bundle $V^*_\Theta$ of dimension
$q-p$ over $\Omega$ with the fiber $V^*_\Theta(w)= (\mbox{ran}\text{
}\Theta(w))^\perp = \mbox{ker} \text{ }\Theta(w)^*$ such that
$$E^*_{\clh_{\Theta}} \cong E^*_{\clh} \otimes V^*_\Theta,$$
where $\clh_{\Theta}=(\clh \otimes \mathbb{C}^q)/ M_{\Theta} (\clh
\otimes \mathbb{C}^p)$. Moreover, $\clh_{\Theta} \in
B^*_{q-p}(\Omega)$. Finally, one has the identity
$$\clk_{E^*_{\clh_{\Theta}}} - \clk_{E^*_{\clh}} \otimes I_{V^*_\Theta} =
I_{E^*_{\clh}} \otimes \clk_{V^*_\Theta}.$$

\end{Theorem}

\NI\textsf{Proof.} Let $\{e_i\}_{i=1}^p$ and $\{\hat{e}_j\}_{j=1}^q$
be the standard orthonormal bases of $\mathbb{C}^p$ and
$\mathbb{C}^q$, respectively. We first show that $\mbox{ker~} (M_z - w
I_{\clh})^* \otimes {\mbox{ker~} \Theta(w)^*} \subseteq \clh \otimes
\mathbb{C}^q$ is orthogonal to $\mbox{ran~} M_{\Theta}$.

\NI For $B = \sum_{i=1}^p b_i \otimes e_i \in \clh \otimes
\mathbb{C}^p$, $k_w \in \mbox{ker~}(M_z - w I_{\clh})^*$ and $\Xi_w =
\sum_{j=1}^q \alpha_j \hat{e}_j \in \mbox{ker~} \Theta(w)^*$, we have
\[\langle M_{\Theta} B, k_w \otimes {\Xi}_w\rangle_{\clh \otimes
\mathbb{C}^q} = \langle B, M_{\Theta}^*(k_w \otimes
{\Xi}_w)\rangle_{\clh \otimes \mathbb{C}^p} = \sum_{j=1}^q
\sum_{i=1}^p b_i(w) \theta_{j,i}(w) {\alpha}_j = 0,\]because
\[\sum_{i=1}^p b_i(w)\{\sum_{j=1}^q \theta_{j,i}(w) {\alpha}_j\} e_i
= {\Theta(w)^* {\Xi}_w} = 0 \in \mathbb{C}^p.\]Next, we show
that \[\bigvee_{w \in \Omega} \{\mbox{ker~}(M_z - w I_{\clh})^* \otimes
{\mbox{ker~}\Theta(w)^*}\} = (\clh \otimes \mathbb{C}^q) \ominus
\mbox{ran~}M_{\Theta}.\]Let $C = \sum_{j=1}^q c_j \otimes \hat{e}_j
\in \clh \otimes \mathbb{C}^q$ be such that $\langle C, k_w \otimes
{\Xi}_w\rangle _{\clh\otimes \mathbb{C}^q} = 0$ for all $k_w \in
\mbox{ker~}(M_z - wI_{\clh})^*$ and $\Xi_w \in
\mbox{ker~}\Theta(w)^*$. We want to show that there exists an $H =
\sum_{i=1}^p h_i \otimes e_i \in \clh \otimes \mathbb{C}^p$ such
that $C = M_{\Theta} H$, which would complete the proof that the
eigenspaces of $M_z$ on $\clh_{\Theta}$ span $\clh_{\Theta}$.

We have $C = M_{\Theta}H$ if and only if $c_j(w) = \sum_{i=1}^p
\theta_{j,i}(w) h_i(w)$ for $j=1, \ldots, p$ and $w \in \Omega$. Let
$w_0 \in \Omega$. Since $\mbox{rank} \{\theta_{j,i}(w_0)\}$ is $p$,
there exist integers $1 \leq j_1 < j_2 \ldots < j_p \leq q$ such that the $p
\times p$ matrix made up of the $p$ rows
$\{\theta_{j_k,i}(w_0)\}_{k=1}^p$ has a non-zero determinant. Hence,
there exists a unique $p$-tuple, denoted by $\{h_i(w_0)\}_{i=1}^p
\in \mathbb{C}^p$, such that $c_j(w_0) = \sum_{i=1}^p
\theta_{j,i}(w_0)h_i(w_0)$ for $j = j_1, \ldots, j_p$. For any other
row $j_0$, there exists a $p$-tuple $\{\gamma_{j_k}\}_{k=1}^p \in \mathbb{C}^p$
such that
\[\theta_{j_0,i}(w_0) = \sum_{k=1}^p \gamma_{j_k}
\theta_{j_k,i}(w_0),\]for $k = 1, \ldots, p$. Thus the vector
$\Gamma = \sum_{k=1}^p {\gamma}_{j_k} \hat{e}_{j_k} - \hat{e}_{j_0}
\in \clh \otimes \mathbb{C}^q$ satisfies $\Theta(w_0)^* {\Gamma} = 0
\in \clh \otimes \mathbb{C}^p$, and hence, $k_{w_0} \otimes {\Gamma} \in
\mbox{ker~} (M_z - w_0 I_{\clh})^* \otimes \mbox{ker~}\Theta(w_0)^*$.
Therefore, \[\langle C, k_{w_0} \otimes {\Gamma}\rangle_{\clh
\otimes \mathbb{C}^q} = 0,\] and hence
\[\begin{split} 0 & = \langle \sum_{j=1}^q c_j \otimes \hat{e}_j,
k_{w_0} \otimes (\sum_{k=1}^p {\gamma}_{j_k} \hat{e}_{j_k} -
\hat{e}_{j_0})\rangle = \sum_{j=1}^p \langle c_{j_k},
k_{w_0}\rangle_{\clh}
\bar{\gamma}_{j_k} - \langle c_{j_0}, k_{w_0}\rangle_{\clh}\\
& = \sum_{k=1}^p c_{j_k}(w_0) \bar{\gamma}_{j_k} -
c_{j_0}(w_0),\end{split}
\]or \[c_{j_0}(w_0) = \sum_{k=1}^p c_{j_k}(w_0)
\bar{\gamma}_{j_k}.\]Thus, \[C(w_0) = \Theta(w_0)
F(w_0).\]Since $w_0 \in \Omega$ is arbitrary, we have defined a
function $H = \sum_{i=1}^p h_i \otimes e_i$ on all of $\Omega$ but
we need to show that $H \in \clh \otimes \mathbb{C}^p$. Recall that
$\Theta$ has a left inverse $\Psi \in \clm(\clh) \otimes
\clb(\mathbb{C}^q, \mathbb{C}^p)$. Thus we have $H(w) =
\Psi(w) \Theta(w) F(w) = \Psi(w) C(w)$ which implies that \[H =
M_{\Psi} C \in \clh \otimes \mathbb{C}^p.\]

The closed range condition for $\clh_{\Theta}$ follows as in the
proof of Theorem \ref{mcase} and hence $\clh_{\Theta} \in
B_1^*(\Omega)$.

To establish the curvature formula, we first recall that the formula
for the curvature of the Chern connection on an open subset $U
\subseteq \Omega$ for a hermitian anti-holomorphic vector bundle is
$\bar{\partial}[G^{-1}{\partial}G]$, where $G$ is the Gramian for an
anti-holomorphic frame $\{f_i\}_{i=1}^{q-p}$ for the vector bundle
on $U$ (cf. \cite{CS}). We assume that $U$ is chosen so that the
$\{k_w\}$ for $w \in \Omega$ can be chosen to be an anti-holomorphic
function on $U$. Denoting by $G_{\Theta}$ the Gramian for the frame
$\{k_w \otimes f_i(w)\}_{i=1}^{q-p}$, $G_{\Theta}(w)$ equals the
$(q-p) \times (q-p)$ matrix
\[G_{\Theta}(w) = \big( \langle k_w \otimes f_i(w), k_w \otimes
f_j(w)\rangle \big)_{i, j = 1}^{q-p} = \|k_w\|^2 \big ( \langle
f_i(w), f_j(w)\rangle \big)_{i,j = 1}^{q-p} = \|k_w\|^2
G_{f}(w),\]where $G_{f}$ is the Gramian for the anti-holomorphic
frame $\{f_i(w)\}_{i=1}^{q-p}$ for $V^*_\Theta$. Then
\[
\begin{split}
\bar{\partial} [G_{\Theta}^{-1}({\partial} G_{\Theta})] &=
\bar{\partial} [\frac{1}{\|k_w\|^2} G_{f}^{-1}({\partial}(\|k_w\|^2
G_{f}))]
\\& =
\bar{\partial} [\frac{1}{\|k_w\|^2} G_{f}^{-1}({\partial}(\|k_w\|^2)
G_{f} + \|k_w\|^2 {\partial} G_{f})] \\& = \bar{\partial}
[\frac{1}{\|k_w\|^2}{\partial}(\|k_w\|^2) + G_{f}^{-1} {\partial}
G_{f} ]\\ & = \bar{\partial}
[\frac{1}{\|k_w\|^2}{\partial}(\|k_w\|^2)] + \bar{\partial} [
G_{f}^{-1}{\partial} G_{f}].
\end{split}
\]
Hence, expressing these matrices in terms of the respective frames
and using the fact that the coordinates of a bundle and of its dual
can be identified using the basis given by the frame, one has
\[\clk_{E^*_{\clh_{\Theta}}}(w) - \clk_{E^*_{\clh}}(w) \otimes
I_{V^*_\Theta(w)} = I_{E^*_\clh (w)} \otimes \clk_{V^*_\Theta}(w),\]
for all $w \in U$. Since the coordinate free formula does not involve
$U$, this completes the proof. \qed

\vspace{0.2in}

Based on the Theorem just stated, we can say that the isomorphism
of quotient Hilbert modules is independent of the choice of the
basic Hilbert module "building blocks" from which they were
created.

\begin{Corollary}\label{n-curv-cor}
Let $\clh, \tilde{\clh} \in B_1^*(\Omega)$ for $\Omega \subseteq
\mathbb{C}^n$. For $i = 1,2$, assume that $\Theta_i : \Omega \raro
\clb(\mathbb{C}^p, \mathbb{C}^q)$ is in both $\clm(\clh) \otimes
\clb(\mathbb{C}^p, \mathbb{C}^q)$ and $\clm(\tilde{\clh})
\otimes \clb(\mathbb{C}^p, \mathbb{C}^q)$. Moreover,
assume that $\Theta_i$ has a left inverse multiplier for both $\clh$
and $\tilde{\clh}$ and for $i=1, 2$. Then $\clh_{\Theta_1}$ is
isomorphic to $\clh_{\Theta_2}$ if and only if
$\tilde{\clh}_{\Theta_1}$ is isomorphic to
$\tilde{\clh}_{\Theta_2}$.
\end{Corollary}

\NI\textsf{Proof.} The statement is obvious from the tensor product
representations $E^*_{\clh_{\Theta_i}} \cong E^*_{\clh} \otimes
V^*_{\Theta_i}$ and $E^*_{\tilde{\clh}_{\Theta_i}} \cong
E^*_{\tilde{\clh}} \otimes V^*_{\Theta_i}$, for $i=1, 2$ ; that is, isomorphic as hermitian anti-holomorphic bundles, and the result that
$\clk_{E^*_{\clh_{\Theta_1}}} = \clk_{E^*_{\clh_{\Theta_2}}}$ if and
only if $\clk_{V_{\Theta_1}} = \clk_{V_{\Theta_2}}$ as two forms. \qed

\vspace{0.1in}

In \cite{DKKS} we showed that for $\clh$ and $\tilde{\clh}$ in the
standard family of contractive Hilbert modules over the disk algebra
$A(\mathbb{D})$ and multipliers $\Theta$ and $\tilde{\Theta}$
that, if $\clh_{\Theta}$ and $\tilde{\clh}_{\tilde{\Theta}}$ are
isomorphic, then so are $\clh$ and $\tilde{\clh}$. (Recall that the
disk algebra $A(\mathbb{D})$ consists of all continuous functions on
the closure of $\mathbb{D}$ that are holomorphic on $\mathbb{D}$ with
the supremum norm.) Therefore, for this family of quotient Hilbert
modules, the isomorphism question reduced to an earlier version of
Corollary \ref{cor-3.2}. It seems possible that such a result might
hold in greater generality. Establishing it, however, would depend
on having a better understanding of how the curvatures of Hilbert
modules in $B_1(\Omega)$ are related to those of the holomorphic
sub-bundles of the product bundles $\Omega \times \mathbb{C}^q$.

\vspace{0.2in}

\newsection{Similarity of Quotient Hilbert modules}

In this section, we investigate conditions for certain quotient
Hilbert modules to be similar to the reproducing kernel Hilbert
modules from which they are constructed. We begin with the case in
which the existence of a left inverse for the multiplier depends
only on a positive answer to the corona problem for the domain.

\begin{Theorem}\label{sim1}
Let $\clh$ be a scalar-valued reproducing kernel Hilbert module over $\mathbb{C} [z_1, \cdots, z_n]$ related to
$\Omega \subseteq \mathbb{C}^n$. Assume that $\theta_1, \theta_2,
\psi_1, \psi_2 $ are in $\clm(\clh)$ and that $\theta_1 \psi_1 +
\theta_2 \psi_2 =1$. Then the quotient Hilbert module
$\clh_{\Theta}= (\clh \otimes \mathbb{C}^2)/ M_{\Theta} \clh$ is
similar to $\clh$, where $M_{\Theta} f = \theta_1 f \otimes e_1 +
\theta_2 f \otimes e_2 \in \clh \otimes \mathbb{C}^2$ and $f \in
\clh$, with $\{e_1, e_2\}$ the standard orthonormal basis for $\mathbb{C}^2$.
\end{Theorem}
\NI\textsf{Proof.} Let $R_{\Psi} : \clh \oplus \clh \raro \clh$ be
the bounded module map defined by $R_{\Psi} (f \oplus g) = \psi_1 f
+ \psi_2 g$ for $f, g \in \clh$. Note that
$$R_{\Psi}M_{\Theta}=I_{\clh},$$ or that $R_{\Psi}$ is a left
inverse for $M_\Theta$. Then for any $f \oplus g \in \clh \oplus
\clh$, we have $$f \oplus g = (M_{\Theta} R_{\Psi}(f \oplus g)) + (
f \oplus g - M_{\Theta} R_{\Psi}(f \oplus g)),$$ with $M_{\Theta}
R_{\Psi}(f \oplus g) \in \mbox{ran}\, M_{\Theta}$ and $f \oplus g -
M_{\Theta} R_{\Psi}(f \oplus g) \in \mbox{ker}\, R_{\Psi}$. This
decomposition, along with $$\mbox{ran}\, M_{\Theta} \cap
\mbox{ker}\, R_{\Psi} = \{0\}$$ implies that $$\clh \oplus \clh =
\mbox{ran}\, M_{\Theta} \stackrel{\cd}+ \mbox{ker}\, R_{\Psi}.$$
Thus, there exists a module idempotent $Q \in \clb(\clh \oplus
\clh)$ with matrix entries in $\clm(\clh)$ such that $Q(\Theta f +
g) =g$ for $f \in \clh$ and $g \in \mbox{ker~}R_{\Psi}$. Moreover,
$\mbox{ran}\, M_{\Theta} = \mbox{ker}\, Q$ and $\mbox{ker}\,
R_{\Psi} = \mbox{ran}\, Q$. The composition $Q \circ \pi_{\Theta}^{-1}
: \clh_{\Theta} \raro \clh$ is well-defined and is the required
invertible module map establishing the similarity of $\clh_\Theta$
and $\clh$.  \qed

It has been observed by earlier authors that the case for $n=2$ is
much simpler than for $n >2$ (cf. \cite{OF}).

\begin{Corollary}\label{sim2}
Let $\theta_1, \theta_2 \in \clm(H^2_n)$ satisfy $|\theta_1(z)|^2 +
|\theta_2(z)|^2 \geq \epsilon$ for all $z \in \mathbb{B}^n$ and some
$\epsilon >0$. Then the quotient Hilbert module
$(H^2_n)_{\Theta}=(H^2_n \otimes \mathbb{C}^2) / M_{\Theta} H^2_n$
is similar to $H^2_n$.
\end{Corollary}
\NI\textsf{Proof.} The corollary follows from Theorem \ref{sim1} using the
corona theorem for $\clm(H^2_n)$ (see \cite{CSW} or \cite{OF}).\qed

\vspace{0.1in}

We next show that the similarity criterion for
a certain class of quotient Hilbert modules is independent of the
choice of the basic Hilbert module ``building blocks" as in the
isomorphism case, so long as the multiplier algebras are the same.

Recall that a short exact sequence of Hilbert modules  \[0
\longrightarrow \clh \otimes \mathbb{C}^p \stackrel{M_\Theta}
{\longrightarrow} \clh \otimes \mathbb{C}^q \stackrel{\pi_{\Theta}}
\longrightarrow \clh_{\Theta} \longrightarrow 0,\] is said to
\emph{split} if $\pi_\Theta$ is right invertible; that is, if there
exists a module map $\sigma_\Theta: \clh_\Theta \raro \clh \otimes
\mathbb{C}^q$ such that
$$\pi_{\Theta}\sigma_{\Theta}=I_{\clh_\Theta}.$$
In the algebraic context, splitting is equivalent to $M_{\Theta}$
being left invertible. In \cite{D2} and \cite{DFS}, this fact was
extended to Hilbert modules with a straightforward proof.

Thus the question of similarity of a quotient Hilbert module to the
building block Hilbert module can be raised in the context of a
split short exact sequence. In general, splitting is not equivalent
to similarity. The question is related to the corona problem (cf.
\cite{D2}) and the commutant lifting theorem (cf. \cite{DFS}).
However, somewhat surprisingly, the relation doesn't depend on the
Hilbert module so long as the multiplier algebra remains the same.

\begin{Theorem}\label{bundle-map}
Let $\clh, \tilde{\clh} \in B_1^*(\Omega)$ for $\Omega \subseteq
\mathbb{C}^n$, be such that $\clm(\clh) \subseteq \clm(\tilde{\clh})$ and let $\Theta \in \clm(\clh) \otimes \clb(\mathbb{C}^p, \mathbb{C}^q)$
, for $1 \leq p <q$, be left invertible. Then the similarity of ${\clh}_{\Theta} = (\clh \otimes \mathbb{C}^q)/
M_{\Theta} (\clh \otimes \mathbb{C}^p)$ to $\clh \otimes \mathbb{C}^{q-p}$ implies the similarity of $\tilde{{\clh}}_{\Theta} = ({\tilde{\clh}} \otimes \mathbb{C}^q)/
M_{\Theta} ({\tilde{\clh}} \otimes \mathbb{C}^p)$ to ${\tilde{\clh}} \otimes \mathbb{C}^{q-p}$.
\end{Theorem}

\NI\textsf{Proof.} Since $\clm(\clh) \subseteq \clm(\tilde{\clh})$,
$\Theta \in \clm(\tilde{\clh}) \otimes \clb(\mathbb{C}^p,
\mathbb{C}^q)$ and $\tilde{\clh}_{\Theta}$ is well-defined.
Moreover, by Theorem \ref{mcase}, we have $\clh_{\Theta},
\tilde{\clh}_{\Theta} \in B^*_{q-p}(\Omega)$.

\NI Let $\sigma_{\Theta}$ be a module cross-section for
$\clh_{\Theta}$; that is, $\sigma_{\Theta}: \clh_{\Theta}
\rightarrow \clh \otimes \mathbb{C}^q$ such that
$$\pi_{\Theta}\sigma_{\Theta}=I_{\clh_{\Theta}}.$$
If we set
$$Q=\sigma_{\Theta}\pi_{\Theta},$$
we obtain a module idempotent on $\clh \otimes \mathbb{C}^q$ such that
$$
\mbox{ran~}Q \stackrel{.}+ \mbox{ran~}M_{\Theta}= \clh \otimes \mathbb{C}^q.
$$
But there exists a $\Phi \in \clm(\clh) \otimes \clb(\mathbb{C}^q)$ such that
$$
M_{\Phi}=Q,
$$
and $\Phi(z)$ is an idempotent on $\clb(\mathbb{C}^q)$ for $z \in \Omega$. An easy argument using localization shows that
$$
\mbox{ran~}\Phi(z) \stackrel{.}+ \mbox{ran~}\Theta(z)= \mathbb{C}^q,
$$
for $z \in \Omega$. But this fact is independent of $\clh$.

\NI Thus, if we set
$$
\tilde{\sigma}_{\Theta}=M_{\Phi}\tilde{\pi}^{-1}_{\Theta}
$$
on $\tilde{\clh} \otimes \mathbb{C}^q$, where $\tilde{\pi}_{\Theta}$ is the quotient map of the short exact sequence for ${\tilde{\clh}}_{\Theta}$, then $\tilde{\sigma}_{\Theta}$ is a module map from $\tilde{\clh}_{\Theta}$ to $\tilde{\clh} \otimes \mathbb{C}^q$. Moreover, the idempotent $\tilde{Q}=\tilde{\sigma}_{\Theta}\tilde{\pi}_{\Theta}$ is again represented by $M_{\Phi}$.

Suppose that $\clh_{\Theta}$ is similar to $\clh \otimes \mathbb{C}^{q-p}$. Then there exists an invertible module map $X: \clh \otimes \mathbb{C}^{q-p} \rightarrow \clh_{\Theta}$. Compose the module maps $\sigma_{\Theta}$ and $X$ to obtain $Y=\sigma_{\Theta}X: \clh \otimes \mathbb{C}^{q-p} \rightarrow \clh \otimes \mathbb{C}^q$ and let $\Gamma \in \clm(\clh) \otimes \clb(\mathbb{C}^{q-p}, \mathbb{C}^q)$ so that $Y=M_{\Gamma}$. Since $\clm(\clh) \subseteq \clm(\tilde{\clh})$,  we can use $\Gamma$ to define
$$
M_{\Gamma}: \tilde{\clh} \otimes \mathbb{C}^{q-p} \rightarrow \tilde{\clh} \otimes \mathbb{C}^q.
$$
Composing $M^{-1}_{\Gamma}$ and $\tilde{\sigma}_{\Theta}$, we obtain an invertible module map $M^{-1}_{\Gamma} \tilde{\sigma}_{\Theta}$ from $\tilde{\clh}_{\Theta}$ to $\tilde{\clh} \otimes \mathbb{C}^{q-p}$, which shows that $\tilde{\clh}_{\Theta}$ is similar to $\tilde{\clh} \otimes \mathbb{C}^{q-p}$. \qed

\begin{Corollary}\label{similar}
Let $\tilde{\clh} \in B^*_1(\mathbb{D})$ be a bounded Hilbert module over $\mathbb{C}[z_1, \cdots, z_n]$; that is, a Hilbert module $\tilde{\clh}$ for which there exists a constant $C$ such that $\|p \cdot f \| \leq C \|p \|_{A(\mathbb{D})}\|h \|_{\tilde{\clh}}$ for all $p \in \mathbb{C}[z_1, \cdots, z_n]$ and $h \in \tilde{\clh}$. Moreover, let  $\Theta \in H^{\infty}(\mathbb{D}) \otimes \clb(\mathbb{C}^p, \mathbb{C}^q)$ have a left inverse in $H^{\infty}(\mathbb{D}) \otimes \clb(\mathbb{C}^q, \mathbb{C}^p)$. Then $\tilde{\clh}_{\Theta}$ is similar to $\tilde{\clh} \otimes \mathbb{C}^{q-p}$.
\end{Corollary}
\NI\textsf{Proof.} We use Theorem \ref{bundle-map} with $\clh=H^2(\mathbb{D})$ and $\tilde{\clh}$ the given Hilbert module. We observe following \cite{D2} that for a bounded Hilbert module related to $\mathbb{D}$, $\clm(\clh)=H^{\infty}(\mathbb{D})$. The proof is completed by appealing to a result of Sz.-Nagy and Foias about a left invertible $\Theta$ (cf. \cite{NF}). \qed

\vspace{0.1in}

The question of similarity is equivalent to a problem in complex geometry (cf. \cite{D2}). In general, for a split short exact sequence
$$0 \longrightarrow \clh \otimes \mathbb{C}^p \stackrel{M_\Theta} {\longrightarrow}
\clh \otimes \mathbb{C}^q \stackrel{\pi_{\Theta}} \longrightarrow
  \clh_{\Theta} \longrightarrow 0,$$
one can define the idempotent function $\Gamma: \Omega \rightarrow \clb(\mathbb{C}^q)$, where $\mbox{ran~}\Gamma$ yields a hermitian holomorphic subbundle $\clf$ of the trivial bundle $\Omega \times \mathbb{C}^q$. If $\Omega$ is contractible, then $\clf$ is trivial. The question of similarity is equivalent to whether one can find a trivializing frame for which the corresponding Gramian $G$ is uniformly bounded above and below when $\clm(\clh)=H^{\infty}(\Omega)$, or  it and its inverse lie in the multiplier algebra when it is smaller. As mentioned above, this question is related to the corona problem and the commutant lifting theorem.

Corollary \ref{similar} along with Theorem 0.1 in \cite{KT}
provide a connection between the quotient Hilbert modules of the
Hardy module and those of any other ``reasonable" reproducing kernel
Hilbert modules over $A(\mathbb{D})$ such as the weighted Bergman
modules. (See \cite{DKT} for a more detailed account of this
phenomenon for the Bergman space based on function theory. The
results there expand on the following result.)

\begin{Corollary}\label{similar-link}
Under the assumptions of Corollary \ref{similar}, the following
statements are equivalent, where $\Pi(z)$ denotes the orthogonal
projection of $\clh_{\Theta}$ onto $\mbox{ker~}(M_z - w
I)^*|_{\clh_{\Theta}}$, the localization of $\pi_{\Theta}$ at $z$. $\clf$ is the same finite-dimensional Hilbert space throughout these statements:

(1) $\clh_\Theta$ is similar to $\clh \otimes \clf$.

(2) $(H^2(\mathbb{D}))_\Theta$ is similar to $H^2(\mathbb{D}) \otimes \clf$.

(3) The eigenvector bundles of $\clh_\Theta$ and $\clh \otimes \clf$
are uniformly equivalent; that is, there exists an anti-holomorphic
pointwise invertible bundle map $\Phi : E^*_{\clh_{\Theta}} \raro
E^*_{\clh} \otimes \clf$ and a scalar $c
> 0$ such that $\frac{1}{c} \|f\| \leq \|\Phi(w)f\| \leq c \|f\|$ for all $w
\in \mathbb{D}$ and $f \in E^*_{\clh_{\Theta}}$.

(4) There exists a bounded subharmonic function $\varphi$ defined on
$\mathbb{D}$ such that
$$
\nabla^2 \vp(z) \geq -{\clk}_{E^*_{\clh_\Theta}}(z)+{\clk}_{E^*_{\clh} \otimes \clf}(z)
$$
for all $z \in \mathbb{D}$.

(5) The measure $$[-{\clk}_{E^*_{\clh_\Theta}}(z)+{\clk}_{E^*_{\clh} \otimes \clf}(z)] (1-|z|) dx dy$$ is a Carleson measure,
and the
estimate
$$-{\clk}_{E^*_{\clh_\Theta}}(z)+{\clk}_{E^*_{\clh} \otimes \clf}(z) \leq \frac{C}{(1-|z|)^2}$$ holds
for some $C >0$.
\end{Corollary}

\NI\textsf{Proof.} One considers the previous Theorem along with the
results in \cite{KT} to establish these equivalences. \qed

\newsection{Concluding Remark}

The results and techniques in this paper raise many questions on
possible extensions and generalizations. We mention a few.

First, one could attempt to generalize Theorem \ref{n-curv} to the
case in which $\Theta(z)$ has constant rank but a nontrivial kernel.
If we assume that the kernel of $\Theta(z)$ has constant dimension,
then two bundles are defined, one by the kernel and the other by the
co-kernel of $\Theta(z)$, and the quotient Hilbert module
$\clh_{\Theta}$ would have a resolution by a tensor product of Hilbert
modules in the form $\clh \otimes \mathbb{C}^l$, for some positive
integers $l$, of longer length depending on resolving the kernel
bundle. Probably in this case, a curvature formula would exist and
involve an alternating sum of curvatures of the modules in the
resolution. These questions are most likely related to the results
in \cite{DFS}.

Second, characteristic operator functions, of which $\Theta(z)$ is
an example, do not, in general, have constant rank. Still the
techniques of this paper might be useful in studying them, at least
for isomorphism questions, since curvatures need to be equal on only
an open set to conclude isomorphism. Similarity would be another
matter, however.

Further, for a Hilbert module $\clh \in B_m^*(\Omega)$, $1 \leq m <
\infty$, and $\Omega \subseteq \mathbb{C}^n$, the multiplier Banach
algebra $\clm(\clh)$ can be identified with the commutant of the
algebra of operators defined by $\mathbb{C}[z_1, \ldots, z_n]$
acting on $\clh$. For many results in this paper concerning two
Hilbert modules, $\clh$ and $\clh_*$, the important assumption is that
$\clm(\clh) = \clm(\clh_*)$. Moreover, if $\clh$ and $\tilde{\clh}$
are in $B^*_m(\Omega)$, the question of whether $\clm(\clh,
\tilde{\clh}) \neq \{0\}$ is closely related to that of the
similarity of $\clh$ and $\tilde{\clh}$, or of ``parts" of them,
which is closely related to the ``similarity" of the corresponding
hermitian anti-holomorphic rank $m$ vector bundles $E_{\clh}^*$ and
$E^*_{\tilde{\clh}}$. A key difficulty in making these relationships
precise concerns the question of when bundle maps between bundles,
obtained as pull backs from the Grassmaninan for a complex Hilbert
space, can be realized as the result of global maps on the Hilbert
space. At present, there are no geometric criteria known
guaranteeing that this is possible. For unitary equivalence, the
situation is much simpler which is one of the key observations in
\cite{CD} allowing one to relate operator theory to complex
geometry.

Finally, a necessary condition for the quotient Hilbert module
$\clh_{\Theta}$ to be similar to $\clh \otimes \mathbb{C}^{q-p}$
(using the notation of the previous section) is that the bundle
$E^*_{\clh_{\Theta}}$ be trivial as a hermitian anti-holomorphic
vector bundle. Therefore, the first step in establishing similarity is
to construct a bounded anti-holomorphic bundle map $\eta : \Omega
\times \mathbb{C}^{q-p} \raro \Omega \times \mathbb{C}^q$ so that
$\mbox{ran~} \eta(z) \stackrel{.}+ \mbox{ran~} \Theta (z) =
\mathbb{C}^q$ for $z \in \Omega$. Next, one would need a module map
$\Xi: \clh \otimes \mathbb{C}^{q-p} \raro \clh \otimes \mathbb{C}^q$
so that $\Xi(z) = \eta(z)$ for $z \in \Omega$.  But to get started
the bundle $E^*_{\clh_{\Theta}}$ must be trivial.

If there exist hermitian anti-holomorphic vector bundles $E$ and $F$
so that $E$ and $E \oplus F$ are trivial with bounded
trivializations but $F$ is not trivial, then we could construct an
example of a quotient Hilbert module for which $E^*_{\clh_{\Theta}}$
is not trivial. This would provide an example for which there is a
topological obstruction to similarity. Can this happen?

For many domains $\Omega$ such as the unit ball or the polydisk, all
vector bundles are trivial.  But there still is a more
subtle possible obstruction to similarity.  In particular, the bundle
could be trivial but not posses a bounded trivialization.  Can this
happen? We mean by a bounded trivialization a bounded bundle map
$\psi: \Omega \times \mathbb{C}^{q-p} \raro E^*_{\clh_{\Theta}}$
with $\frac{1}{c}\|f\| \leq \|\Psi(z)f\| \leq c\|f\|$ for some
$c>0$, and for all $f \in \mathbb{C}^{q-p}$ and $z \in \Omega$.

A detailed discussion of the relation of these notions is in
\cite{D2} and the case $\clh = H^2_n$ is given at the beginning of
Section 3 in \cite{DFS}.

\end{document}